\input xy \xyoption{all}

\documentclass[12pt,leqno]{article}

\usepackage{amsfonts}
\usepackage{amsmath}
\usepackage{amssymb}

\title{\sc Homotopy types of orbit spaces and their self-equivalences
for the  periodic groups ${\mathbb Z}/a \rtimes ({\mathbb
Z}/b\times T^\star_n)$ and ${\mathbb Z}/a\rtimes({\mathbb
Z}/b\times O^\star_n)$}

\date      \kill

\author{Marek Golasi\'nski and Daciberg Lima Gon\c{c}alves}

\renewcommand{\k}{\hspace{-2mm}{\bf .}\hspace{5mm}}

\newtheorem{thmX}{Theorem}[section]
\newtheorem{proX}[thmX]{Proposition}

\newtheorem{CorX}[thmX]{Corollary}

\renewcommand{\baselinestretch}{1.2}

\begin{document}

\maketitle

\vspace*{-7mm}






{\bf Abstract.} {\footnotesize Let $G$ be a finite group given in
one of the forms listed in the title with period $2d$ and $X(n)$
an $n$-dimensional $CW$-complex with the homotopy type of an
$n$-sphere.
\par We study the automorphism group $\mbox{Aut}\,(G)$ to compute the number
of distinct homotopy types of orbit spaces $X(2dn-1)/\mu$ with
respect to free and cellular $G$-actions $\mu$ on all
$CW$-complexes $X(2dn-1)$. At the end, the groups ${\mathcal
E}(X(2dn-1)/\mu)$ of self homotopy equivalences of orbit spaces
$X(2dn-1)/\mu$ associated with free and cellular $G$-actions $\mu$
on $X(2dn-1)$ are determined.}

\renewcommand{\baselinestretch}{1.33}

\vspace{10mm}

{\large\bf Introduction.} Given a free and cellular action $\mu$ of a finite group $G$ with order $|G|$
on a $CW$-complex $X$, write $X/\mu$ for the corresponding
orbit space. The problem of determining all possible homotopy types
of $X/\mu$ among all free and cellular actions $\mu$ on $X$,
as well the group ${\mathcal E}(X/\mu)$ of self homotopy equivalences
of $X/\mu$ has been extensively studied for a number of spaces e.g., in \cite{R}.
Notoriously, for an odd dimensional sphere $\mathbb{S}^{2n-1}$ with a free action
of a finite cyclic group $\mathbb{Z}/k$ this corresponds to the classification
of lens spaces and the calculation of the groups of it self homotopy
equivalences studied in \cite{GG}. A larger family of interesting examples are given by a free and cellular action of a finite group $G$ with order $|G|$
on a $CW$-complex $X(2n-1)$ with the homotopy type of a $(2n-1)$-sphere. Write $X(2n-1)/\mu$ for the corresponding
orbit space called a $(2n-1)$-{\em spherical space form} or a {\em Swan $(2n-1)$-complex} (see e.g., \cite{AD}).
Taking into account \cite{Sm}, the case of spherical space forms presents a special interest. Furthermore,
Swan \cite{Sw} has shown that any finite group with periodic cohomology of period $2d$ acts freely
and cellularly on a $(2d-1)$-dimensional $CW$-complex of the homotopy type of a $(2d-1)$-sphere. It is worth to mention that useful cohomological and geometric aspects associated to group actions
are presented in \cite{AD} and a list of basic conjectures is provided.
\par Backing to the case of a $(2n-1)$-dimensional $CW$-complex $X(2n-1)$
with the homotopy type of a $(2n-1)$-sphere, by means of results in \cite{Sw}, it is shown in \cite[Theorem 1.8]{Th1} that the set of homotopy types of spherical
space forms of all free cellular $G$-actions on $X(2n-1)$ is in one-to-one correspondence with the orbits, which contain a generator
of the cyclic group $H^{2n}(G)={\mathbb Z}/{|G|}$ under the action of $\pm\mbox{Aut}\,(G)$ (see \cite{GG} for another approach). This plays also
a fundamental role in the calculation of the group $\mathcal{E}(X(2n-1)/\mu)$ of self homotopy equivalences
of the orbit space $X(2n-1)/\mu$.
\par All finite periodic groups has been completely described by Suzuki-Zassenhaus and their classification can be found in
the table \cite[Chapter IV; Theorem 6.15]{AM}. The present paper is part of the project to describe the homotopy
types of the orbit spaces and the group of self homotopy equivalences for all periodic groups. It continues the works
of \cite{GG,GG1,GG2,GG4}, where the cases corresponding to the families I
and II from the table \cite[Chapter IV; Theorem 6.15]{AM} with the Suzuki-Zassenhaus
classification of finite periodic groups have been solved. Here we have two goals.
The first one is to calculate the numbers of homotopy types of spherical spaces forms
for the groups $\mathbb{Z}/a\rtimes(\mathbb{Z}\times T^\star_n)$ and
$\mathbb{Z}/a\rtimes(\mathbb{Z}\times O^\star_n)$ corresponding to the families III
and IV from the table mentioned above. The second one is to determine the group of homotopy
classes of self-equivalences for space forms given by free actions of those both families of finite
periodic groups. The results of \cite{GG,GG1,GG2,GG4}, taking care for the groups from families
I and II of that table, are essential to make crucial calculations to develop the main results
stated in Theorem 2.2 and Theorem 3.2.
\par In order to obtain these results, we divide the paper into two parts.
The first part consists of some algebraic results. The automorphism group $\mbox{Aut}(A\rtimes_\alpha G)$ of a semi-direct product $A\rtimes_\alpha G$
of some finite groups $A,G$ leads in \cite{GG2} to a splitting short exact sequence
$$0 \to \mbox{Der}_\alpha\,(G,A)\longrightarrow\mbox{Aut}\,(A\rtimes_\alpha G)\longrightarrow\mbox{Aut}\,(A)\times\mbox{Aut}_\alpha(G)\to 1.$$
Section 1 makes use of this to achieve automorphisms of the groups in question.
This is the approach to develop in Proposition 1.1 and Proposition 1.2 the groups
$\mbox{Der}_\alpha\,(G,A)$ and $\mbox{Aut}\,(A)\times\mbox{Aut}_\alpha(G)$, respectively.
\par Then, in the second part, we present geometric interpretations of those algebraic results in terms of $G$-actions.
Section 2 uses the group $\mbox{Aut}\,(A)\times\mbox{Aut}_\alpha(G)$ established in Proposition 1.2 and Lyndon-Hochschild-Serre spectral sequence to deal with the number of homotopy
types of spherical space forms for actions of the groups
${\mathbb Z}/a\rtimes({\mathbb Z}/b\times T^\star_n)$ and ${\mathbb Z}/a\rtimes
({\mathbb Z}/b\times O^\star_n)$. The main results of this section are stated in

{\bf Theorem 2.2.} {\em  Let $\gamma=(\gamma_1,\gamma_2) : {\mathbb Z}/b\times T^\star_n\to({\mathbb Z}/a)^\star$ and
$\tau=(\tau_1,\tau_2) : {\mathbb Z}/b\times O^\star_n\to({\mathbb Z}/a)^\star$ be actions
with $(a,b)=(ab,6)=1$ and $n\ge 3$, where $\gamma_1 : {\mathbb Z}/b\to({\mathbb Z}/a)^\star$,
$\gamma_2 : T^\star_n\to({\mathbb Z}/a)^\star$ and $\tau_1: \mathbb{Z}/b\to(\mathbb{Z}/a)^\star$,
$\tau_2 : O^\star_n\to(\mathbb{Z}/a)^\star$ are appropriate restrictions
of $\gamma$ and $\tau$, respectively. Then$:$

\vspace{3mm}

{\em (1)} $\mbox{\em card}\,{\mathcal K}^{2k[\ell(\gamma),2]-1}_{\mathbb{Z}/a\rtimes_\gamma
(\mathbb{Z}/b\times T^\star_n)}/_\simeq=2^{t+t'+1}3^{n_0}O(a,k[\ell(\gamma),2])
O_{\mbox{\tiny\em Aut}_{\gamma_1}(\mathbb{Z}/b)}(b,k[\ell(\gamma),2])O(3^{n-n_0},$

$k[\ell(\gamma),2])^{-1}$ for some $0\le t\le 2$ and $0\le t'\le 1$$;$

\vspace{3mm}

{\em (2)} $\mbox{\em card}\,{\mathcal K}^{2k[\ell(\tau),2]-1}_{\mathbb{Z}/a\rtimes_\tau(\mathbb{Z}/b\times O^\star_n)}/_\simeq=2^{t+1}\times 3^{n-1}
O(a,k[\ell(\tau),2])O_{\mbox{\tiny\em Aut}_{\tau_1}(\mathbb{Z}/b)}(b,k[\ell(\tau),2])$ for

some $0\le t\le 1$.}

Then, Corollary 2.3 says that the number of such
homotopy types of those space forms coincides with that of $(4n-1)$-lens spaces
studied in \cite{GG} provided the least period of the groups in question is $\le 4$.

\par The group of crossed homomorphisms $\mbox{Der}_\alpha\,(G,A)$ studied in Proposition 1.1
plays a key role in Section 3 dealing with the structure of groups ${\mathcal E}(X(2dn-1)/\mu)$
of self homotopy equivalences for spherical space forms $X(2dn-1)/\mu$ with respect to free and cellular
${\mathbb Z}/a\rtimes({\mathbb Z}/b\times T^\star_n)$-- and
${\mathbb Z}/a\rtimes({\mathbb Z}/b\times O^\star_n)$--actions $\mu$, respectively.
We point out that by means of \cite[Proposition 3.1]{GG} (see also \cite[Theorem 1.4]{Sm}), the group
${\mathcal E}(X(2k-1)/\mu)$ is independent of the action $\mu$ on $X(2k-1)$. Writing $X(2k-1)/G$ for
the corresponding orbit space, we close the paper with

{\bf Theorem 3.2.}  {\em Let $\gamma=(\gamma_1,\gamma_2) : \mathbb{Z}/b\times T^\star_n\to(\mathbb{Z}/a)^\star$ $($resp.\
$\tau=(\tau_1,\tau_2) : \mathbb{Z}/b\times O^\star_n\to(\mathbb{Z}/a)^\star$$)$ be an action with
$(a,b)=(ab,6)=1$ for $n\ge 3$.
If the group ${\mathbb Z}/a\rtimes_\gamma({\mathbb Z}/b\times
T^\star_n)$ $($resp.\ ${\mathbb Z}/a\rtimes_\tau({\mathbb Z}/b\times O^\star_n)$$)$
acts freely and cellularly on a $CW$-complex $X(2k[\ell(\gamma),2]-1)$
$($resp.\ $X(2k[\ell(\tau),2]-1)$$)$ then
\begin{multline*}{\mathcal E}(X(2k[\ell(\gamma),2]-1)/({\mathbb Z}/a\rtimes_\gamma({\mathbb Z}/b\times
T^\star_n)))\cong\mbox{\em Der}_\gamma\,
({\mathbb Z}/b\times T^\star_n,{\mathbb Z}/a)\rtimes\\({\mathcal E}(X(2k[\ell(\gamma_1),2]-1)/({\mathbb Z}/a))\times
{\mathcal E}_{\gamma_1}\,(X(2k[\ell(\beta),2]-1)/({\mathbb Z}/b))\times
S_4\times\\\mathbb{Z}/\bigg(\frac{3^{n-n_0}}{(3^{n-n_0},k[\ell(\gamma),2])}\bigg)\\\end{multline*}

\vspace{-2cm}

\begin{multline*}(\mbox{resp.}\;\;{\mathcal E}(X(2k[\ell(\tau),2]-1)/({\mathbb Z}/a\rtimes_\tau({\mathbb Z}/b\times O^\star_n)))\cong\mbox{\em Der}_\tau\,
({\mathbb Z}/b\times O^\star_n,{\mathbb Z}/a)\rtimes\\({\mathcal E}(X(2k[\ell(\tau_1),2]-1)/({\mathbb Z}/a))\times
{\mathcal E}_{\tau_1}\,(X(2k[\ell(\beta),2]-1)/({\mathbb Z}/b))\times
O_n\rtimes\\\mathbb{Z}/\bigg(\frac{3^{n-1}}{(3^{n-1},k[\ell(\tau),2])}\bigg))\\\end{multline*}}

\vspace{-1.7cm}

\noindent
which deals with explicit formulae for those groups of self homotopy equivalences.
\par Approaching of homotopy types of spherical space forms and their self homotopy equivalences
for the rest of the groups from the table in \cite[Chapter IV; Theorem 6.15]{AM}, or more precisely for
the family of VI of this table, is in progress.
\par {\small The authors are grateful to the referee for carefully reading earlier version of
the paper and all his suggestions to make the introduction clear and understandable.
The main part of this work has been done during the visit of the first author to
the Department of Mathematics-IME, University of S\~ao Paulo during the period July 09--August 08, 2003.
He would like to thank the Department of Mathematics-IME for its hospitality during his stay. This visit was
supported by FAPESP, Projecto Tem\'atico Topologia Alg\'ebrica, Geom\'etrica e Differencial-2000/05385-8, Ccint-USP
and Projecto 1-Pr\'o-Reitoria de Pesquisa-USP.}

 \newpage

\renewcommand{\baselinestretch}{1.27}

\large

\normalsize

\renewcommand{\thesection}{{}}
\section{}
\renewcommand{\thesection}{\arabic{section}}

{\large\bf 1.\ Algebraic backgrounds.} Let a finite group $G$ be given by an
extension
$$1 \to G_1 \to G \to G_2\to 1,$$
where the orders of groups $G_1$ and $G_2$ are relatively prime. We recall that
by \cite{GG3} any automorphism of $G$ leaves the subgroup $G_1$ invariant and consequently, there is a map
$\psi :\mbox{Aut}\,(G)\to\mbox{Aut}\,(G_1)\times\mbox{Aut}\,(G_2)$ of automorphism groups.
Given an $H$-action $\alpha : H\to \mbox{Aut}\,(A)$ on an abelian group $A$
write $\mbox{Der}_\alpha\,(H,A)$ for the abelian group of crossed homomorphisms.
For $H$-actions $\alpha_1 : H\to\mbox{Aut}\,(A_1)$ and
$\alpha_2 : H\to\mbox{Aut}\,(A_2)$ consider the obvious induced action $(\alpha_1,\alpha_2) : H\to\mbox{Aut}\,(A_1\times A_2)$.
Then, an isomorphism $$ \mbox{Der}_{(\alpha_1,\alpha_2)}\,(H,A_1\times A_2)\stackrel{\cong}{\longrightarrow}
\mbox{Der}_{\alpha_1}\,(H,A_1)\times\mbox{Der}_{\alpha_2}\,(H,A_2)\leqno(\star)$$
follows.
\par Now, let $0\to A \to G \to H \to 1$ be a short exact sequence, with $A$ an
abelian group. Then, there is an obvious $H$-action $\alpha : H\to \mbox{Aut}\,(A)$.
If groups $A$ and $H$ are finite with relatively prime orders then the cohomology
group $H^1(H,A)$ vanishes (see e.g., \cite[Corollary 5.4]{AM}) and consequently, $\mbox{Der}_\alpha\,(H,A)=A/A^H$,
where $A^H$ is the subgroup of $A$ consisting of all elements fixed  under the action of $H$.
Furthermore, by \cite[Lemma 1.2]{GG2} this sequence $0 \to A \to G\to H \to 1$
of finite groups yields the exact sequence
$$0 \to \mbox{Der}_\alpha\,(H,A)\to\mbox{Aut}(G)\stackrel{\psi}{\to}\mbox{Aut}\,(A)\times\mbox{Aut}\,(H).$$
\par For an action $\alpha : G\to\mbox{Aut}\,(A)$, let $A\rtimes_\alpha G$
denote the semi-direct product of $A$ and $G$ with respect to the action
$\alpha$. Let the orders of $A$ and $G$ be relatively primes  and
$\psi:\mbox{Aut}(A\rtimes_\alpha G)
\to \mbox{Aut}\,(A)\times\mbox{Aut}(G)$ be the obvious map. Then, by \cite{GG2},
$\mbox{Im}\,\psi=\mbox{Aut}\,(A)\times\mbox{Aut}_\alpha(G)$, where
$\varphi\in\mbox{Aut}_\alpha(G)$ if and only if $\alpha=\alpha\varphi$ or
equivalently  $$\mbox{Aut}_\alpha(G)=\{\varphi\in\mbox{Aut}\,(G);\;
\varphi(\mbox{Ker}\alpha)=\mbox{Ker}\,\alpha\;\mbox{and}\;
\bar{\varphi}=\mbox{id}_{G/\mbox{\tiny Ker}\,\alpha}\},$$ where $\bar{\varphi}$
denotes the map induced by $\varphi$ on the quotient group $G/\mbox{Ker}\,\alpha$.

\bigskip

\par Now, let $Q_8$ be the classical quaternion group $\{\pm 1,\pm i,\pm j,\pm k\}$
of order $8$, where $1,i,j$ and $k$ are generators of the quaternion algebra over
reals. Consider the action $\alpha : \mathbb{Z}/3\to \mbox{Aut}\,(Q_8)$ such that
a generator of $\mathbb{Z}/3$ is sent to the automorphism $\tau\in \mbox{Aut}\,(Q_8)$
defined by: $\tau(i)=j$, $\tau(j)=k$ and $\tau(k)=i$. Since $\mbox{Aut}\,(Q_8){\cong}
S_4$, the symmetric group on four letters (see e.g., \cite[Lemma 6.9]{AM}) any two
faithful representations of $\mathbb{Z}/3$ in the group $Q_8$ are conjugated.
Whence, without losing generality, we can choose the action $\alpha$ given
above. Then, we consider the  semi-direct product $Q_8\rtimes_\alpha \mathbb{Z}/3=T^\star$,
the {\em binary tetrahedral} group. More generally, for $n\ge 1$ consider the
action $\alpha_n : \mathbb{Z}/{3^n}\to \mbox{Aut}\,(Q_8)$ as the composition of the quotient map $\mathbb{Z}/3^n\to\mathbb{Z}/3$ with
the action $\alpha : \mathbb{Z}/3\to \mbox{Aut}\,(Q_8)$.
Then, for the group $$T^\star_n=Q_8\rtimes_{\alpha_n}\mathbb{Z}/{3^n},$$
by means of \cite[p.\ 198]{Wo}, it holds
$$T^\star_n:\begin{cases}X^{3^n}=P^4=1,\, P^2=Q^2,\, XPX^{-1}=Q,\cr
XQX^{-1}=PQ,\,PQP^{-1}=Q^{-1}
\end{cases}$$ in virtue of generators and relations.
In particular, the cyclic group $\mathbb{Z}/3^n$ is the abelianization of $T^\star_n$ for any $n\ge 1$
and the center $\mathcal{Z}(T^\star_n)=\mathbb{Z}/2\oplus\mathbb{Z}/3^{n-1}$.
\par  The symmetric group $S_3$ has two distinct extensions by $Q_8$, with respect to
the outer action $\alpha : S_3 \to \mbox{Out}\,(Q_8)=\mbox{Aut}\,(Q_8)/\mbox{Inn}\,(Q_8)$ which is
the composition of the inclusion $S_3\subseteq \mbox{Aut}\,(Q_8)$ with the projection
$\mbox{Aut}\,(Q_8) \to \mbox{Out}(Q_8)$. This follows from the facts that $\mathcal{Z}(Q_8)=\mathbb{Z}/2$ and $H^2(S_3,\mathbb{Z}/2)=\mathbb{Z}/2$.
These extensions are the  semi-direct product
$Q_8\rtimes S_3$ and $$1\to Q_8\to O^\star\stackrel{\varphi}{\to}S_3\to 1,$$
where $O^\star$ is the {\em binary octahedral} group. Because $\varphi^{-1}(A_3)=T^\star$
for the alternating subgroup $A_3\subseteq S_3$, so we achieve the extension
$$1\to T^\star\to O^\star\to \mathbb{Z}/2\to 1.$$
In general, since $\mathcal{Z}(T^\star_n)=\mathbb{Z}/2\oplus\mathbb{Z}/3^{n-1}$
and $H^2(S_3,\mathcal{Z}(T^\star_n))=H^2(S_3,\mathbb{Z}/2\oplus\mathbb{Z}/3^{n-1})= H^2(S_3,\mathbb{Z}/2)=\mathbb{Z}/2$,
we achieve the non-trivial extension $$1\to T^\star_{n-1}\to O^\star_n\stackrel{\varphi_n}{\to}S_3\to 1$$ for $n\ge 1$, where $T^\star_0=Q_8$. Because $\varphi_n^{-1}(A_3)=T^\star_n$ a {\em fortiori} the new extension $$1\to T^\star_n\to O^\star_n\to\mathbb{Z}/2\to 1$$ is obtained.
\par In the light of \cite[p.\ 198]{Wo} the group $O^\star_n$ is given by
$$O^\star_n:\begin{cases}X^{3^n}=P^4=1,\,P^2=Q^2=R^2,\,PQP^{-1}=Q^{-1},\cr
XPX^{-1}=Q,\,XQX^{-1}=PQ,\,RXR^{-1}=X^{-1},\\
RPR^{-1}=QP,\,RQR^{-1}=Q^{-1}
\end{cases}$$ in virtue of generators and relations.
It follows that the cyclic group $\mathbb{Z}/2$ is isomorphic to the abelianization of $O^\star_n$ and the
center $\mathcal{Z}(O^\star_n)$ as well.
\par Now, consider the periodic groups ${\mathbb Z}/a\rtimes_\gamma({\mathbb Z}/b\times T^\star_n)$ and
${\mathbb Z}/a\rtimes_\tau({\mathbb Z}/b\times O^\star_n)$ corresponding to the families
III and IV \cite[Theorem 6.15]{AM} with $(a,b)=(ab,6)=1$ and $n\ge 1$, where
$\gamma : {\mathbb Z}/b\times T^\star_n\to \mbox{Aut}\,({\mathbb Z}/a)$ and
$\tau : {\mathbb Z}/b\times O^\star_n\to \mbox{Aut}\,({\mathbb Z}/a)$ are actions of ${\mathbb Z}/b\times T^\star_n$
and  ${\mathbb Z}/b\times O^\star_n$, respectively, on the cyclic group ${\mathbb Z}/a$.
The group $\mbox{Aut}\,({\mathbb Z}/a)$ is abelian, a {\em fortiori} the actions
$\gamma$ and $\tau$ are uniquely determined by their restrictions
$\gamma_1 : {\mathbb Z}/b\to \mbox{Aut}\,({\mathbb Z}/a)$, $\gamma_2 :
T^\star_n\to \mbox{Aut}\,({\mathbb Z}/a)$ and $\tau_1 : {\mathbb Z}/b\to \mbox{Aut}\,({\mathbb Z}/a)$,
$\tau_2 : O^\star_n\to \mbox{Aut}\,({\mathbb Z}/a)$. But the abelianizations of
$T^\star_n$ and $O^\star_n$ are isomorphic to the groups ${\mathbb Z}/3^n$ and
$\mathbb{Z}/2$, respectively. Whence, the actions $\gamma_2$ and $\tau_2$ are
uniquely determined by $\gamma_2(X)$ and $\tau_2(R)$, respectively, with
$\gamma_2(X)^{3^n}=\mbox{id}_{\mathbb{Z}/a}$ and $\tau_2(R)^2=\mbox{id}_{\mathbb{Z}/a}$.

\par To study the groups $\mbox{Aut}\,({\mathbb Z}/a\rtimes_\gamma({\mathbb Z}/b\times T^\star_n))$ and $\mbox{Aut}\,({\mathbb Z}/a\rtimes_\gamma({\mathbb Z}/b\times O^\star_n))$,
we need, in the light of \cite[Proposition 1.3]{GG2}, to describe the groups
$\mbox{Der}_\gamma\,({\mathbb Z}/b\times T^\star_n,{\mathbb Z}/a)$, $\mbox{Aut}_\gamma\,({\mathbb Z}/b\times T^\star_n)$ and
$\mbox{Der}_\tau\,({\mathbb Z}/b\times O^\star_n,{\mathbb Z}/a)$, $\mbox{Aut}_\tau\,({\mathbb Z}/b\times O^\star_n)$,
respectively.
\par First, let $a=p^k$ with  $p\not=2,3$ prime and $k\ge 0$. Because the actions
$\gamma_2$ and $\tau_2$ factor through the abelianizations of $T^\star_n$ and
$O^\star_n$, which are isomorphic to the groups $\mathbb{Z}/3^n$ and $\mathbb{Z}/2$, respectively,
whence $\mbox{Ker}\,\gamma_2$ is trivial or $\mbox{Ker}\,\gamma_2=Q_8\rtimes_{\alpha_n}\mathbb{Z}/3^{n_0}$ for some $n_0\le n$ and
$\mbox{Ker}\,\tau_2$ is trivial or equals $T^\star_n$ (as a subgroup with index two in
$O^\star_n$ and containing $T^\star_n$). Consequently, by \cite{GG2}, we achieve that $\mbox{Der}_\gamma\,({\mathbb Z}/b\times T^\star_n,{\mathbb Z}/p^n)
=\mbox{Der}_{\gamma_1}\,({\mathbb Z}/b,{\mathbb Z}/p^n)$ provided $\gamma_2$ is trivial and $\mbox{Der}_\gamma\,({\mathbb Z}/b\times T^\star_n,{\mathbb Z}/p^k)=\mbox{Der}_{\bar{\gamma}}\,({\mathbb Z}/b\times{\mathbb Z}/3^{n-n_0},{\mathbb Z}/p^n)$ provided $\mbox{Ker}\,\gamma_2=Q_8\rtimes_{\alpha_n}\mathbb{Z}/3^{n_0}$, where $\bar{\gamma} : {\mathbb Z}/b\times{\mathbb Z}/3^{n-n_0}={\mathbb Z}/b\times (T^\star_n/\mbox{Ker}\,\gamma_2)\to\mbox{Aut}\,({\mathbb Z}/a)$ is the action induced by $\gamma$.
Furthermore, $\mbox{Der}_\tau\,({\mathbb Z}/b\times O^\star_n,{\mathbb Z}/p^n)
=\mbox{Der}_{\tau_1}\,({\mathbb Z}/b,{\mathbb Z}/p^n)$ provided $\tau_2$ is trivial and $\mbox{Der}_\tau\,({\mathbb Z}/b\times O^\star_n,{\mathbb Z}/p^k)=\mbox{Der}_{\bar{\tau}}\,({\mathbb Z}/b\times{\mathbb Z}/2,{\mathbb Z}/p^n)$ provided $\mbox{Ker}\,\tau_2=T^\star_n$,
where $\bar{\tau} : {\mathbb Z}/b\times{\mathbb Z}/2={\mathbb Z}/b\times O^\star_n/T^\star_n\to\mbox{Aut}\,({\mathbb Z}/a)$ is the action induced by $\tau$.
Because $(b,6)=1$, the groups ${\mathbb Z}/b\times{\mathbb Z}/3^{n-n_0}$ and ${\mathbb Z}/b\times{\mathbb Z}/2$ are cyclic whence, as in \cite[Corollary 1.5]{GG2}, elements of $\mbox{Der}_{\bar{\tau}}\,({\mathbb Z}/b\times{\mathbb Z}/3^{n-n_0},{\mathbb Z}/p^n)$
and $\mbox{Der}_{\bar{\tau}}\,({\mathbb Z}/b\times{\mathbb Z}/2,{\mathbb Z}/p^n)$ might be described by means of some elements in ${\mathbb Z}/p^n$.
\par  Now, if $a$ is a positive integer with $(a,6)=1$ and $a=p_1^{k_1}\cdots p_s^{k_s}$ its prime factorization with $k_i\ge 1$ then $p_i\not=2,3$ for all $i=1,\ldots,s$. Obviously, any isomorphism ${\mathbb Z}/a\stackrel{\cong}{\to}{\mathbb Z}/{p_1^{n_1}}\times\cdots\times{\mathbb Z}/{p_s^{n_s}}$
yields isomorphisms $\alpha :\mbox{Aut}\,({\mathbb Z}/a)\stackrel{\cong}{\to}\mbox{Aut}\,({\mathbb Z}/{p_1^{n_1}}\times\cdots\times{\mathbb Z}/{p_s^{n_s}})$ and
$$\mbox{Der}_\gamma\,({\mathbb Z}/b\times T^\star_n,{\mathbb Z}/a)\stackrel{\cong}{\to}\mbox{Der}_{\alpha\gamma}({\mathbb Z}/b\times
T^\star_n,{\mathbb Z}/{p_1^{k_1}}\times\cdots\times{\mathbb Z}/{p_s^{k_s}})$$
for an action $\gamma : {\mathbb Z}/b\times T^\star_n\to\mbox{Aut}\,({\mathbb Z}/a)$, and
$$\mbox{Der}_\tau\,({\mathbb Z}/b\times O^\star_n,{\mathbb Z}/a)\stackrel{\cong}{\to}\mbox{Der}_{\alpha\tau}({\mathbb Z}/b\times
O^\star_n,{\mathbb Z}/{p_1^{k_1}}\times\cdots\times{\mathbb Z}/{p_s^{k_s}})$$
for an action $\tau : {\mathbb Z}/b\times O^\star_n\to\mbox{Aut}\,({\mathbb Z}/a)$.
Then, the well-known (see e.g.\ \cite[Lemma 1.1]{GG2}) isomorphism $\mbox{Aut}\,({\mathbb Z}/{p_1^{k_1}}\times\cdots\times{\mathbb Z}/{p_s^{k_s}})\stackrel{\cong}{\to}
\mbox{Aut}\,({\mathbb Z}/{p_1^{k_1}})\times\cdots\times\mbox{Aut}\,({\mathbb Z}/{p_s^{k_s}})$ and $(\star)$
lead to isomorphisms
$$\mbox{Der}_\gamma\,({\mathbb Z}/b \times T^\star_n,{\mathbb Z}/a)\stackrel{\cong}{\longrightarrow}
\mbox{Der}_{\alpha_1\gamma}({\mathbb Z}/b\times T^\star_n,{\mathbb Z}/{p_1^{k_1}})\times\cdots\times\mbox{Der}_{\alpha_s\gamma}\,
(\mathbb{Z}/b\times T^\star_n,{\mathbb Z}/{p_s^{n_s}})$$
and
$$\mbox{Der}_\tau\,({\mathbb Z}/b \times O^\star_n,{\mathbb Z}/a)\stackrel{\cong}{\longrightarrow}
\mbox{Der}_{\alpha_1\tau}({\mathbb Z}/b\times O^\star_n,{\mathbb Z}/{p_1^{k_1}})\times\cdots\times\mbox{Der}_{\alpha_s\tau}\,
(\mathbb{Z}/b\times O^\star_n,{\mathbb Z}/{p_s^{n_s}}),$$
where $\alpha_i$ is the composition of $\alpha$ with an appropriate projection
map $\mbox{Aut}\,({\mathbb Z}/{p_1^{k_1}})\times\cdots\times\mbox{Aut}\,({\mathbb Z}/{p_s^{k_s}})\to\mbox{Aut}\,({\mathbb Z}/p_i^{k_i})$ for $i=1,\ldots,s$.
Thus, we may summarize the discussion above as follows.
\begin{proX}\k Let ${\mathbb Z}/b$ and ${\mathbb Z}/p^k$ be cyclic groups with
 $p$ prime and $k\ge 1$,
$(bp^k,6)=(b,p^k)=1$ and let $\gamma : {\mathbb Z}/b\times T^\star_n\to\mbox{\em Aut}\,({\mathbb Z}/p^k)$,
$\tau : {\mathbb Z}/b\times O^\star_n\to\mbox{\em Aut}\,({\mathbb Z}/p^k)$ be actions.
Write $\gamma_1 : \mathbb{Z}/b\to\mbox{\em Aut}\,({\mathbb Z}/p^k)$,
$\gamma_2 : T^\star_n\to\mbox{\em Aut}\,({\mathbb Z}/p^k)$
and $\tau_1 : \mathbb{Z}/b\to\mbox{\em Aut}\,({\mathbb Z}/p^k)$, $\tau_2 : O^\star_n\to\mbox{\em Aut}\,({\mathbb Z}/p^k)$ for the appropriate restrictions of $\gamma$ and $\tau$, respectively.
Then:

\vspace{2mm}

{\em (1)} $$\mbox{\em Der}_\gamma({\mathbb Z}/b\times T^\star_n,{\mathbb Z}/p^k)\cong\mbox{\em Der}_{\gamma_1}\,
({\mathbb Z}/b,{\mathbb Z}/p^k)$$ and $$\mbox{\em Der}_\tau({\mathbb Z}/b\times O^\star_n,{\mathbb Z}/p^k)
\cong\mbox{\em Der}_{\tau_1}\,({\mathbb Z}/b,{\mathbb Z}/p^k)$$
if $\gamma_2$ and $\tau_2$ are trivial;

\vspace{2mm}

{\em (2)} $\mbox{\em Der}_\gamma({\mathbb Z}/b\times T^\star_n,{\mathbb Z}/p^k)\cong\mbox{\em Der}_{\bar{\gamma}}\,({\mathbb Z}/b\times{\mathbb Z}/3^{n-n_0},{\mathbb Z}/p^k)$
provided $\mbox{\em Ker}\,\gamma_2=Q_8\rtimes_{\alpha_n}\mathbb{Z}/3^{n_0}$, where
$\bar{\gamma} : {\mathbb Z}/b\times{\mathbb Z}/3^{n-n_0}\cong{\mathbb Z}/b\times(T^\star_n/\mbox{\em Ker}\,\gamma_2)\to
\mbox{\em Aut}\,({\mathbb Z}/p^k)$ is the action induced by $\gamma$\newline

\noindent
and

$\mbox{\em Der}_\tau({\mathbb Z}/b \times O^\star_n,{\mathbb Z}/p^k)\cong\mbox{\em Der}_{\bar{\tau}}\,({\mathbb Z}/b\times{\mathbb Z}/2,{\mathbb Z}/p^k)$
provided $\mbox{\em Ker}\,\tau_2=T^\star_n$, where
$\bar{\tau} : {\mathbb Z}/b\times{\mathbb Z}/2\cong{\mathbb Z}/b\times(O_n^\star/T^\star_n)\to
\mbox{\em Aut}\,({\mathbb Z}/p^k)$ is the action induced by $\tau$. \label{D}

\vspace{2mm}

\par If $\gamma : {\mathbb Z}/b\times T^\star_n\to \mbox{\em Aut}\,({\mathbb Z}/a)$ and  $\tau : {\mathbb Z}/b\times O^\star_n\to \mbox{\em Aut}\,({\mathbb Z}/a)$ are actions
with $(a,b)=(ab,6)=1$ and $a=p_1^{k_1}\cdots p_s^{k_s}$ is the prime factorization of $a$ with
$k_i\ge 1$ for $i=1,\ldots,s$ then
$$\mbox{\em Der}_\gamma\,({\mathbb Z}/b\times T^\star_n,{\mathbb Z}/a)\stackrel{\cong}{\longrightarrow}
\mbox{\em Der}_{\alpha_1\gamma}({\mathbb Z}/b\times T^\star_n,{\mathbb Z}/{p_1^{k_1}})\times\cdots\times
\mbox{\em Der}_{\alpha_s\gamma}\,(\mathbb{Z}/b\times T^\star_n,\mathbb{Z}/{p_s^{n_s}})$$ and
$$\mbox{\em Der}_\tau\,({\mathbb Z}/b\times O^\star_n,{\mathbb Z}/a)\stackrel{\cong}{\longrightarrow}
\mbox{\em Der}_{\alpha_1\tau}({\mathbb Z}/b\times O^\star_n,{\mathbb Z}/{p_1^{k_1}})\times\cdots\times
\mbox{\em Der}_{\alpha_s\tau}\,({\mathbb Z}/b\times T^\star_n,{\mathbb Z}/{p_s^{n_s}}),$$ where
$\alpha_i$ is the composition of an isomorphism  $\alpha: \mbox{\em Aut}\,({\mathbb Z}/a)
\stackrel{\cong}{\to}\mbox{\em Aut}\,({\mathbb Z}/{p_1^{k_1}}\times\cdots\times{\mathbb Z}/{p_s^{n_s}})$
with an appropriate projection map $\mbox{\em Aut}\,({\mathbb Z}/{p_1^{k_1}})\times\cdots\times
\mbox{\em Aut}\,({\mathbb Z}/{p_s^{n_s}})\to\mbox{\em Aut}\,({\mathbb Z}/p_i^{k_i})$ for $i=1,\ldots,s$.
\end{proX}

\par Now, move to the groups $\mbox{Aut}_\gamma({\mathbb Z}/b\times T^\star_n)$ and $\mbox{Aut}_\tau({\mathbb Z}/b\times O^\star_n)$, where $\gamma=(\gamma_1,\gamma_2) :
{\mathbb Z}/b\times T^\star_n\to\mbox{Aut}\,({\mathbb Z}/a)$ and $\tau=(\tau_1,\tau_2) :
{\mathbb Z}/b\times O^\star_n\to\mbox{Aut}\,({\mathbb Z}/a)$.
 Because $(ab,6)=1$, \cite[Lemma 1.1]{GG2} yields $\mbox{Aut}\,({\mathbb Z}/b\times T^\star_n)\cong
\mbox{Aut}\,({\mathbb Z}/b)\times \mbox{Aut}\,(T^\star_n)$ and $\mbox{Aut}\,({\mathbb Z}/b\times O^\star_n)\cong
\mbox{Aut}\,({\mathbb Z}/b)\times \mbox{Aut}\,(O^\star_n)$. Furthermore, the groups $\mbox{Aut}\,(T^\star_n)$ and $\mbox{Aut}\,(O^\star_n)$ have
been fully described in \cite{GG3} for all $n\ge 1$. In the light of \cite[Corollary 1.4]{GG2} we achieve
isomorphisms $$\mbox{Aut}_\gamma\,({\mathbb Z}/b\times T^\star_n)\stackrel{\cong}{\longrightarrow}
\mbox{Aut}_{\gamma_1}\,({\mathbb Z}/b)\times\mbox{Aut}_{\gamma_2}\,(T^\star_n)$$ and
 $$\mbox{Aut}_\tau\,({\mathbb Z}/b\times O^\star_n)\stackrel{\cong}{\longrightarrow}
\mbox{Aut}_{\tau_1}\,({\mathbb Z}/b)\times\mbox{Aut}_{\tau_2}\,(O^\star_n).$$
But $\varphi\in\mbox{Aut}_{\gamma_2}(T^\star_n)$ (resp.\ $\varphi\in\mbox{Aut}_{\tau_2}(O^\star_n)$)
if and only if $\gamma_2(X)=(\gamma_2\varphi)(X)$ (resp.\  $\tau_2(R)=(\tau_2\varphi)(R)$).
Now, if $\mbox{Ker}\,\gamma_2=Q_8\rtimes_{\alpha_n}\mathbb{Z}/3^{n_0}$ then, from the list of elements in $\mbox{Aut}\,(T^\star_n)$
presented in \cite{GG3}, it follows that  $$\mbox{Aut}_{\gamma_2}(T^\star_n)=\{\varphi\in\mbox{Aut}\,(T^\star_n);\;
\varphi(X)=X^{l(1+3^{n_0+1})}\;\mbox{for}\;l=0,\ldots,3^{n-n_0-1}\}.$$
By means of \cite{GG3}, any $\varphi\in\mbox{Aut}\,(O^\star_n)$ restricts to
an automorphism of $T^\star_n$ with the identity on the quotient
$O^\star_n/T^\star_n=\mathbb{Z}/2$ a {\em fortiori} $\tau_2(R)=(\tau_2\varphi)(R)$ holds
for all $\varphi\in\mbox{Aut}\,(O^\star_n)$. Now, in virtue of \cite[Proposition 3.2]{GG3}, we are ready to close this section with

\begin{proX}\k Let ${\mathbb Z}/a$ and ${\mathbb Z}/b$ with $(a,b)=(ab,6)=1$ and $\gamma : {\mathbb Z}/b\times T^\star_n\to
\mbox{\em Aut}\,({\mathbb Z}/a)$,  $\tau : {\mathbb Z}/b\times O^\star_n\to
\mbox{\em Aut}\,({\mathbb Z}/a)$ be actions. Write \label{AQ} $\gamma_1 : {\mathbb Z}/b\to\mbox{\em Aut}\,({\mathbb Z}/a)$,
$\gamma_2 : T^\star_n\to\mbox{\em Aut}\,({\mathbb Z}/a)$ for the restrictions
of $\gamma$ and  $\tau_1 : {\mathbb Z}/b\to\mbox{\em Aut}\,({\mathbb Z}/a)$,
$\tau_2 : O^\star_n\to\mbox{\em Aut}\,({\mathbb Z}/a)$
for the restrictions of $\tau$. Then$:$

\vspace{2mm}

{\em (1)} $\mbox{\em Aut}_{\gamma_2}\,(T^\star_n)\cong S_4\times\mathbb{Z}/3^{n-n_0}$ provided $\mbox{\em Ker}\,\gamma_2=Q_8\rtimes_{\alpha_n}\mathbb{Z}/3^{n_0}$;

\vspace{2mm}

{\em (2)} $\mbox{\em Aut}_{\tau_2}\,(O^\star_n)\cong \mbox{\em Aut}\,(O^\star_n)$.
\end{proX}
Certainly, the groups $\mbox{Aut}_{\gamma_1}\,({\mathbb Z}/b)$ and  $\mbox{Aut}_{\tau_1}\,({\mathbb Z}/b)$ could be described by \cite[Proposition 1.5 and Corollary 1.6]{GG2}. Observe that $\ell(\gamma_1),\ell(\tau_1)\le 2$ implies $\ell(\gamma_1),\ell(\tau_1)=1$ because $b$ is odd and consequently,
$\mbox{Aut}_{\gamma_1}\,({\mathbb Z}/b)=\mbox{Aut}_{\tau_1}\,({\mathbb Z}/b)=\mbox{Aut}\,({\mathbb Z}/b)$, where $\ell(\gamma_1)$ (resp.\ $\ell(\tau_1)$) denotes the order of $\gamma_1(1_b)$ (resp.\ $\tau_1(1_b)$) in $\mbox{Aut}\,(\mathbb{Z}/b)$  for a generator $1_b$ of the cyclic group $\mathbb{Z}/b$.

\renewcommand{\thesection}{{}}
\section{}
\renewcommand{\thesection}{\arabic{section}}

{\large\bf 2.\ Homotopy types of space forms.} Given a group $G$, write $H^k(G)$
 for its $k$th cohomology group with constant coefficients in the integers
${\mathbb Z}$ for $k\ge 0$. Then, any automorphism
$\varphi\in\mbox{Aut}\,(G)$ yields the induced automorphism
$\varphi^\ast\in\mbox{Aut}\, (H^n(G))$ and we write
$\eta : \mbox{Aut}\,(G)\to\mbox{Aut}\,(H^k(G))$ for the corresponding
anti-homomorphism. By a {\em period} of a group $G$ we mean an integer $d$ such
that $H^k(G)=H^{k+d}(G)$ for all $k > 0$, and a group
$G$ with this property is called {\em periodic}. Among all periods of a group
$G$ there is the least one; and all
others are multiple of that one. That least one period we call {\em the period}
of the group  and by \cite[Section 11]{CE} the period of any periodic group is even.
\par Throughout the rest of the paper, $X(k)$ denotes a $k$-dimensional
$CW$-complex with the homotopy type of a $k$-sphere and the group
$\mbox{Aut}\,({\mathbb Z}/a)$ is identified with the unit group
$({\mathbb Z}/a)^\star$ of the $\bmod\, a$ ring ${\mathbb Z}/a$.
Given a free cellular action $\mu$ of a finite group $G$ with order $|G|$
on a $CW$-complex $X(2k-1)$ write $X(2k-1)/\mu$ for the corresponding orbit
space called a $(2k-1)$-{\em spherical space form} or a
{\em Swan $(2k-1)$-complex} (see e.g., \cite{AD}).
Then, the group $G$ is periodic with period $2d$ dividing $2k$ and by
\cite[Chap.\ XVI, \S 9]{CE} there
is an isomorphism $H^{2n}(G)\cong{\mathbb Z}/{|G|}$.
Two spherical space forms $X(2k-1)/\mu$ and $X'(2k-1)/\mu'$ are called
equivalent if they are homeomorphic and let ${\mathcal K}_G^{2k-1}$ denote the
set of all such classes.
We say that two such classes $[X(2k-1)/\mu]$ and $[X'(2k-1)/\mu']$ are
homotopic if the space forms $X(2k-1)/\mu$ and $X'(2k-1)/\mu'$
are homotopy equivalent. Write ${\mathcal K}_G^{2k-1}/_\simeq$ for the associated
quotient set of ${\mathcal K}_G^{2k-1}$
and $\mbox{card}\,{\mathcal K}_G^{2k-1}/_\simeq$ for its cardinality, respectively.
By means of \cite{Sw}, it is shown in \cite[Theorem 1.8]{Th1} that elements of the set
${\mathcal K}_G^{2k-1}/_\simeq$ are in one-to-one correspondence with the orbits,
which contain a generator of $H^{2k}(G)={\mathbb Z}/{|G|}$ under an action of
$\pm\mbox{Aut}\,(G)$
(see also \cite{GG} for another approach). But generators of the group ${\mathbb Z}/{|G|}$ are
given by the unit group $({\mathbb Z}/{|G|})^\star$ of the ring ${\mathbb Z}/{|G|}$. Thus, those
homotopy types are in one-to-one correspondence with the quotient $({\mathbb Z}/{|G|})^\star/\{\pm\varphi^\ast ;\: \varphi\in\mbox{Aut}\,(G)\}$,
where $\varphi^\ast$ is the induced automorphism on the cohomology $H^{2k}(G)={\mathbb Z}/{|G|}$
by $\varphi\in\mbox{Aut}\,(G)$.

\par Now, let $G_1$ and $G_2$ be finite groups with relatively prime orders $|G_1|$
and $|G_2|$, respectively.
If $G_1$ and $G_2$ are also periodic with periods $2d_1$ and $2d_2$, respectively then by \cite{AM} the least
common multiple $[2d_1,2d_2]$ of $2d_1$ and $2d_2$ is the least period of the product $G_1\times G_2$.
Furthermore, given a finite group $G$ with an action $\alpha : G\to(\mathbb{Z}/a)^\star$
write $|\alpha(g)|$ for the order of $\alpha(g)$ with $g\in G$.
Let $\ell(\alpha)=[|\alpha(g)|;\;\mbox{for}\;g\in G]$ be the least common multiple of those orders. Then, for a semi-direct product
$\mathbb{Z}/a\rtimes_\alpha G$, we have shown in \cite{GG2}, by means of the Lyndon-Hochschild-Serre spectral sequence, the following result.
\begin{proX}\k Let $\mathbb{Z}/a$ be a cyclic group of order $a$, $G$ a finite
group, $\alpha : G\to\mathbb{Z}/a$ an action and $(|G|,a)=1$.
If $G$ is periodic with the period $2d$ then the semi-direct product
$\mathbb{Z}/a\rtimes_\alpha G$ is also a period finite group with the least
period $2[\ell(\alpha),d]$.\label{P1}
\end{proX}

\bigskip

\par Thus, we are in a position to investigate the periodic groups ${\mathbb Z}/a\rtimes_\gamma({\mathbb Z}/b\times T^\star_n)$
and ${\mathbb Z}/a\rtimes_\tau({\mathbb Z}/b\times O^\star_n)$.
First, we find the least periods of the groups $T^\star_n$ and $O^\star_n$. Because $T^\star_n=Q_8\rtimes_\alpha\mathbb{Z}/3^n$
whence the Lyndon-Hochschild-Serre spectral sequence applied to the short one
$$0\to Q_8\longrightarrow T^\star_n\longrightarrow \mathbb{Z}/3^n\to 0$$
yields
$$E^{p,q}_2(T_n^\star)=H^p(\mathbb{Z}/3^n,H^q(Q_8))=\begin{cases}0,\;\mbox{if}\; p,q>0;\\
\mathbb{Z},\;\mbox{if}\,p,q=0;\\
0,\;\mbox{if}\;q=0\;\mbox{and}\;p\;\mbox{odd};\\
\mathbb{Z}/3^n,\;\mbox{if}\;q=0\,\mbox{and}\;p\;\mbox{even}\;\mbox{with}\;\not=0;\\
H^0(\mathbb{Z}/3^n,H^q(Q_8))=(H^q(Q_8))^{\mathbb{Z}/3^n},\;\mbox{if}\;q>0.
\end{cases}$$
Using the cohomology
$$H^k(Q_8)=\begin{cases}\mathbb{Z},\;k=0;\\
0,\;\mbox{if}\;k=1+4l;\\
\mathbb{Z}/2\oplus\mathbb{Z}/2,\;\mbox{if}\;k=2+4l;\\
0,\;\mbox{if}\;k=3+4l;\\
\mathbb{Z}/8,\;\mbox{if}\;k=4+4l\end{cases}$$
with $l\ge 0$, we can easily get
$$H^k(T^\star_n)=\begin{cases}\mathbb{Z},\;\mbox{if}\;k=0;\\
0,\;\mbox{if}\;k=1+4l;\\
\mathbb{Z}/3^n,\;\mbox{if}\;k=2+4l;\\
0,\;\mbox{if}\;k=3+4l;\\
\mathbb{Z}/(8\times 3^n),\;\mbox{if}\;k=4+4l
\end{cases}$$
with $l\ge 0$ and consequently, $4$ is the least period of the group $T^\star_n$.
Whence, by Proposition \ref{P1}, the number $2[\ell(\gamma),2]$ is the least period of the
group ${\mathbb Z}/a\rtimes_\gamma({\mathbb Z}/b\times T^\star_n)$.

\bigskip

\par To find the least period of the group $O^\star_n$, we apply Lyndon-Hochschild-Serre spectral sequence to the short one
$$0\to T^\star_n\longrightarrow O^\star_n\longrightarrow \mathbb{Z}/2\to 0.$$
Then, $E^{p,q}_2(O^\star_n)=H^p(\mathbb{Z}/2,H^q(T^\star_n))$. Next, observe that $E^{p,4}_2=H^p(\mathbb{Z}/2,\mathbb{Z}/(8\times 3^n))=H^p(\mathbb{Z}/2,\mathbb{Z}/8)\oplus H^p(\mathbb{Z}/2,\mathbb{Z}/3^n)$. Because $H^p(\mathbb{Z}/2,\mathbb{Z}/3^n)=0$ for $p>0$
and by \cite{Sw} the action of $\mathbb{Z}/2$ on $\mathbb{Z}/8$ is trivial $H^p(\mathbb{Z}/2,\mathbb{Z}/(8\times 3^n))=\mathbb{Z}/2$ for $p>0$. Then, we can easily find that
$$E^{p,q}_2(O^\star_n)=H^p(\mathbb{Z}/2,H^q(T_n^\star))=\begin{cases}\mathbb{Z},\;\mbox{if}\;p=q=0;\\
0,\;\mbox{if}\;p\;\mbox{odd},\, q=0;\\
\mathbb{Z}/2,\;\mbox{if}\;p\;\mbox{even},\, q=0;\\
0,\;\mbox{if}\;p\ge 0,\,q=1+4l,2+4l,3+4l;\\
\mathbb{Z}/(8\times 3^n),\;\mbox{if}\;p=0,\,q=4+4l;\\
\mathbb{Z}/2,\;\mbox{if}\;p>0,\,q=4+4l\end{cases}$$
\noindent with $l \ge 0$.
\par To find the cohomology $H^\ast(O^\star_n)$ consider the generalized quaternion
group $Q_{16}$ as the subgroup of $O^\star_n$ generated by $P,Q,R$ (according to
the presentation of $O^\star_n$ given in Section 1) and its subgroup $Q_8$
generated by $P,Q$. The exact sequence
$0\to Q_8\longrightarrow Q_{16}\longrightarrow\mathbb{Z}/2\to 0$ leads to
Lyndon-Hochschild-Serre spectral sequence with
$E^{p,q}_2(Q_{16})=H^p(\mathbb{Z}/2,H^q(Q_8))$. Because the action of
$\mathbb{Z}/2$ on $\mathbb{Z}/2\oplus\mathbb{Z}/2$ is given by the matrix
$\begin{pmatrix}1&0\\1&1\end{pmatrix}$
and by means of \cite{Sw}, the action of $\mathbb{Z}/2$ on $\mathbb{Z}/8$ is
trivial, applying $H^\star(Q_8)$, we derive:
$$E^{p,0}_2(Q_{16})=H^p(\mathbb{Z}/2)=\begin{cases}\mathbb{Z},\;\mbox{if}\;p=0;\\
0,\;\mbox{if}\;p\;\mbox{odd};\\
\mathbb{Z}/2,\;\mbox{if}\;p\;\mbox{even},\end{cases}$$
$$E^{p,1}_2(Q_{16})=0,\;\;E^{p,2}_2(Q_{16})=H^p(\mathbb{Z}/2,\mathbb{Z}/2\oplus\mathbb{Z}/2)=\begin{cases}\mathbb{Z}/2,\;\mbox{if}\;p=0;\\
0,\;\mbox{if}\; p>0,\end{cases}$$ $E^{p,3}_2(Q_{16})=0$ and $E^{p,4}_2(Q_{16})=H^p(\mathbb{Z}/2,\mathbb{Z}/8)=\mathbb{Z}/2$.
Writing $E_k(Q_{16})$ for the $k$-term of that spectral sequence, we can deduce that
$E_2(Q_{16})\cong E_3(Q_{16})\cong E_4(Q_{16})\cong E_5(Q_{16})$ and $d_5(E^{1,4}_5(Q_{16}))=E^{6,0}_5(Q_{16})$, $d_k(E^{0,q}_k(Q_{16}))=0$
for $k\ge 2$. Then, using the multiplicative structure of that spectral sequence
and the periodicity of the groups $Q_8$ and $\mathbb{Z}/2$, we get further isomorphisms
$E_6(Q_{16})=H(E_5(Q_{16}),d_5)\cong E_7(Q_{16})\cong\cdots\cong E_\infty(Q_{16})
=\mathcal{G}\,(H^\ast(Q_{16}))$, where by \cite[Chapter XII]{CE} it holds:
$$H^k(Q_{16})=\begin{cases}\mathbb{Z},\;\mbox{if}\;k=0;\\
0,\;\mbox{if}\;k=1+4l;\\
\mathbb{Z}\oplus\mathbb{Z}/2,\;\mbox{if}\;k=2+4l;\\
0,\;\mbox{if}\;k=3+4l;\\
\mathbb{Z}/16,\;\mbox{if}\;k=4+4l\end{cases}$$
with $l\ge 0$. The commutative diagram
$$\xymatrix{0\ar[r]&Q_8\ar[r]\ar[d]&Q_{16}\ar[d]\ar[r]&\mathbb{Z}/2\ar@{=}[d]\ar[r]&0\\
0\ar[r]&T^\star_n\ar[r]&O^\star_n\ar[r]&\mathbb{Z}/2\ar[r]&0}$$
leads to a map $$E_k(O^\star_n)\longrightarrow E_k(Q_{16})$$ for $k\ge 2$.
Because of the isomorphism  $E_2^{p,q}(O^\star_n)\stackrel{\cong}{\longrightarrow} E^{p,q}_2(Q_{16})$
for $p>0$, we can get $E_6(O^\star_n)$ and then $E_\infty(O^\star_n)$ as well.
Therefore, we can read that
$$H^k(O^\star_n)=\begin{cases}\mathbb{Z},\;\mbox{if}\;k=0;\\
0,\;\mbox{if}\;k=1+4l;\\
\mathbb{Z}/2,\;\mbox{if}\;k=2+4l;\\
0,\;\mbox{if}\;k=3+4l;\\
A\oplus\mathbb{Z}/3^n,\;\mbox{if}\;k=4+4l\end{cases}$$ with $l\ge 0$, where $A$
is an abelian group of order $16$. Because of the monomorphism
$H^{4l}(O^\star_n)_{(2)}\to H^{4l}(Q_{16})$ on the $2$-primary component of
$H^{4l}(O_n^\star)$ for $l>0$, we deduce an isomorphism $A\cong\mathbb{Z}/16$. Thus,
$$H^k(O^\star_n)=\begin{cases}\mathbb{Z},\;\mbox{if}\;k=0;\\
0,\;\mbox{if}\;k=1+4l;\\
\mathbb{Z}/2,\;\mbox{if}\;k=2+4l;\\
0,\;\mbox{if}\;k=3+4l;\\
\mathbb{Z}/(16\times 3^n),\;\mbox{if}\;k=4+4l\end{cases}$$
with $l\ge 0$ and consequently, $4$ is the least period of the group $O^\star_n$.
Whence, by means of Proposition \ref{P1}, the number $2[\ell(\tau),2]$ is the
least period of the group ${\mathbb Z}/a\rtimes_\tau({\mathbb Z}/b\times O^\star_n)$.

\bigskip

\par  By \cite[Lemma 1.1]{GG2} any automorphism $\varphi\in\mbox{Aut}\,({\mathbb Z}/a\rtimes_\alpha G)$ for $(a,|G|)=1$ determines
a pair $(\varphi_1,\varphi_2)\in({\mathbb Z}/a)^\star\times\mbox{Aut}\,(G)$ with the commutative
diagram
 $$\xymatrix{0\ar[r]&{\mathbb Z}/a\ar[r]\ar[d]^{\varphi_1}&{\mathbb Z}/a\rtimes_\alpha G\ar[d]^\varphi\ar[r]&G\ar[d]^{\varphi_2}\ar[r]&0\\
0\ar[r]&{\mathbb Z}/a\ar[r]&{\mathbb Z}/a\rtimes_\alpha G\ar[r]&G\ar[r]&0.}$$
Then, Lyndon-Hochshild-Serre spectral sequence and its naturality lead to
the commutative diagram of cyclic groups with exact and splitting  rows
$$\xymatrix{0\ar[r]&H^{2k[d,\ell(\alpha)]}(G)\ar[r]\ar[d]^{\varphi^\ast_2}&
H^{2k[d,\ell(\alpha)]}({\mathbb Z}/a\rtimes_\alpha G)\ar[d]^{\varphi^\ast}\ar[r]& H^{2k[d,\ell(\alpha)]}({\mathbb Z}/a)\ar[d]^{\varphi^\ast_1}\ar[r]&0\\
0\ar[r]&H^{2k[d,\ell(\alpha)]}(G)\ar[r]&H^{2k[d,\ell(\alpha)]}({\mathbb Z}/a\rtimes_\alpha G)
\ar[r]&H^{2k[d,\ell(\alpha)]}({\mathbb Z}/a)\ar[r]&0}$$
for $k >0$, where $2d$ is the least period of $G$. Whence, $\varphi^\ast$ is uniquely determined by the corresponding pair $(\varphi^\ast_2,\varphi^\ast_1)$
and consequently, there is the factorization
$$\xymatrix{\mbox{Aut}\,({\mathbb Z}/a\rtimes_\alpha G)\ar[dd]_{\psi}\ar[rr]^\eta&&\mbox{Aut}\,(H^{2k[\ell(\alpha),d]}({\mathbb Z}/a\rtimes_\alpha G))\\
\\
({\mathbb Z}/a)^\star\times\mbox{Aut}_\alpha\,(G)\ar[uurr]_{\eta'}}
$$ for all $k > 0$, where $\mbox{Aut}_\alpha\,(G)$ is the subgroup of $\mbox{Aut}(G)$
defined in Section 1. But $H^{2k[\ell(\alpha),d]}({\mathbb Z}/a\rtimes_\alpha G)\cong{\mathbb Z}/a|G|$,
so in the light of the above, to describe the number $\mbox{card}\,{\mathcal K}_{{\mathbb Z}/a\rtimes_\alpha G}^{2k[\ell(\alpha),d]-1}/_\simeq$
of homotopy types of spherical space forms for ${\mathbb Z}/a\rtimes_\alpha G$ we are led to compute the order
of the quotient $({\mathbb Z}/a|G|)^\star/\{\pm\varphi^\ast;\;\varphi\in({\mathbb Z}/a)^\star\times\mbox{Aut}_\alpha\,(G)\}$
where $\varphi^\ast$ is the induced automorphism on the cohomology $H^{2k[\ell(\alpha),d]}({\mathbb Z}/a\rtimes G)$
for $\varphi\in({\mathbb Z}/a)^\star\times\mbox{Aut}_\alpha\,(G)$.
\par Now, for a periodic group $G_1$ with the least period $2d_1$ and an action
$\omega : G_2\to\mbox{Aut}\,(G_1)$, we achieve anti-homomorphism $G_2\to\mbox{Aut}\,(H^{2kd_1}(G_1))=\mbox{Aut}\,(\mathbb{Z}/|G_1|)$.
Write $(\mathbb{Z}/|G_1|)^\star/\pm G_2$ for the quotient group $(\mathbb{Z}/|G_1|)^\star/\{\pm \omega(g_2)^\ast;\;g_2\in G_2\}$
and $O_{G_2}(|G_1|,2kd_1)$ for its order, where $\omega(g_2)^\ast$ denotes the induced
map on the cohomology $H^{2kd_1}(G_1)$. Furthermore, we set $O(m,n)$ for the order
of the quotient group $(\mathbb{Z}/m)^\star/\{\pm l^n;\;l\in(\mathbb{Z}/m)^\star\}$.
\par Given $\varphi\in\mbox{Aut}\,(T^\star_n)$ for the group
$T^\star_n=Q_8\rtimes_\alpha\mathbb{Z}/3^n$ there is
the corresponding pair
$(\varphi_1,\varphi_2)\in\mbox{Aut}\,(Q_8)\times(\mathbb{Z}/3^n)^\star$ and
by means of \cite{GG3} maps $\varphi_2$ exhaust all automorphisms of the
group $\mathbb{Z}/3^n$. The periodicity of $T^\star_n$,
Lyndon-Hochshild-Serre spectral sequence and
its naturality lead to the commutative diagram
$$\xymatrix{0\ar[r]&H^{4k}(\mathbb{Z}/3^n)\ar[r]\ar[d]^{\varphi^\ast_2}&
H^{4k}(T^\star_n)\ar[d]^{\varphi^\ast}\ar[r]& H^{4k}(Q_8)\ar[d]^{\varphi^\ast_1}\ar[r]&0\\
0\ar[r]&H^{4k}(\mathbb{Z}/3^n)\ar[r]&H^{4n}(T^\star_n)
\ar[r]&H^{4k}(Q_8)\ar[r]&0}$$
of cyclic groups with exact rows for $k\ge 0$.
But, by means of \cite{Sw}, $\varphi_1^\ast$ is the identity map and
a {\em fortiori} $\varphi^\ast$ is uniquely determined by $\varphi_2^\ast$.
\par By \cite{GG3}, any $\varphi\in\mbox{Aut}\,(O^\star_n)$ yields also a pair
$(\varphi_1,\varphi_2)\in\mbox{Aut}\,(T^\star_n)\times(\mathbb{Z}/2)^\star$.
Again, the periodicity of $O^\star_n$, Lyndon-Hochshild-Serre spectral
sequence and its naturality lead to the commutative diagram of cyclic groups
$$\;\xymatrix{0\ar[r]&H^{4k}(\mathbb{Z}/2)\ar[r]\ar[d]^{\varphi^\ast_2}&
H^{4k}(O^\star_n)\ar[d]^{\varphi^\ast}\ar[r]&
H^{4k}(T^\star_n)\ar[d]^{\varphi^\ast_1}\ar[r]&0\\
0\ar[r]&H^{4k}(\mathbb{Z}/2)\ar[r]&H^{4k}(O^\star_n)
\ar[r]&H^{4k}(T^\star_n)\ar[r]&0}$$
with exact and splitting rows for $k\ge 0$. Because the restriction of
$ \varphi_1^\ast$ to $Q_8$, denoted by $\varphi_|$
induces  the identity on $H^{4n}(Q_8)=\mathbb{Z}/8$, by means of the description
of $H^{4k}(O^\star_n)$ and $H^{4k}(T^\star_n)$, we derive from the above the
commutative diagram
$$\;\xymatrix{0\ar[r]&\mathbb{Z}/2\ar[r]\ar@{=}[d]&
\mathbb{Z}/16\ar[d]^{\varphi^\ast_|}\ar[r]&\mathbb{Z}/8\ar@{=}[d]\ar[r]&0\\
0\ar[r]&\mathbb{Z}/2\ar[r]&\mathbb{Z}/16
\ar[r]&\mathbb{Z}/8\ar[r]&0.}$$ Therefore, the restriction
$\varphi^\ast_| : \mathbb{Z}/16\to\mathbb{Z}/16$ is the identity map or
the multipication by $9$. Certainly, both cases might hold. Namely, consider the automorphism $\varphi : O^\star_n\to O^\star_n$
given by $\varphi(P)=P$, $\varphi(Q)=Q$, $\varphi(R)=-R$ and $\varphi(X)=X$,
where $P,Q,R,X$ are generators of $O^\star_n$. But the subgroup of $O^\star_n$
generated by $P,Q,R$ is the generalized quaternion group $Q_{16}$ with the
relation $(RP)^4=R^2$.
Then, $\varphi(RP)=(RP)^5$ and by \cite{Sw} we achieve that $\varphi^\ast : H^{4k}(O^\star_n)\to H^{4k}(O^\star_n)$ restricts on
$H^{4n}(Q_{16})=\mathbb{Z}/16$ to the multiplication by $9$.

\par Now, the group $\mathbb{Z}/a\rtimes_\gamma(\mathbb{Z}/b\times T^\star_n)$ with
an action $\gamma=(\gamma_1,\gamma_2) : \mathbb{Z}/b\times T^\star_n\to(\mathbb{Z}/a)^\star$
yields the factorization
$$\xymatrix{\mbox{Aut}\,({\mathbb Z}/a\rtimes_\gamma(\mathbb{Z}/b\times T^\star_n)) \ar[dd]_{\psi}\ar[rr]^\eta&&\mbox{Aut}\,(H^{2k[\ell(\gamma),2]}({\mathbb Z}/a\rtimes_\gamma(\mathbb{Z}/b\times T^\star_n)) \\
\\
({\mathbb Z}/a)^\star\times\mbox{Aut}_{\gamma_1}(\mathbb{Z}/b)\times\mbox{Aut}_{\gamma_2}(T^\star_n)\ar[uurr]_{\eta'}}
$$ for all $k > 0$. By \cite[Proposition 2.2]{GG2} we obtain the short exact
sequence
$$0\to (\mathbb{Z}/2)^t\longrightarrow (\mathbb{Z}/(8\times 3^nab)/
\pm((\mathbb{Z}/a)^\star\times\mbox{Aut}_{\gamma_1}\,
(\mathbb{Z}/b)\times\mbox{Aut}_{\gamma_2}(T^\star_n))\longrightarrow$$
$$((\mathbb{Z}/a)^\star/\pm(\mathbb{Z}/a)^\star)\times((\mathbb{Z}/b)^\star/\pm\mbox{Aut}_{\gamma_1}(\mathbb{Z}/b))\times((\mathbb{Z}/
(8\times 3^n))^\star/\pm\mbox{Aut}_{\gamma_2}(T^\star_n))\to 0$$
for some $0\le t\le 2$.
But, by means of Proposition \ref{AQ}, $\mbox{Aut}_{\gamma_2}(T^\star_n)\cong S_4\times\mathbb{Z}/3^{n-n_0}$
with $\mbox{Ker}\,\gamma_2=Q_8\rtimes\mathbb{Z}/3^{n-n_0}$. Because the action of $S_4$ on $H^{4k}(Q_8)$
is trivial, the canonical imbedding $(\mathbb{Z}/3^{n-n_0})^\star\hookrightarrow(\mathbb{Z}/3^n)^\star$
leads to the other short exact sequence
$$0\to\mathbb{Z}/2)^{t'}\to(\mathbb{Z}/(8\times 3^n))^\star/
\pm\mbox{Aut}_{\gamma_2}(T^\star_n)\longrightarrow((\mathbb{Z}/8)^\star/\{\pm 1\})
\times((\mathbb{Z}/3^n)^\star/(\mathbb{Z}/3^{n-n_0})^\star)\to 0$$
for some $0\le t'\le 1$.

\bigskip

\par Now, we  move to the group $\mathbb{Z}/a\rtimes_\tau(\mathbb{Z}/b\times O^\star_n)$
with an action $\tau : \mathbb{Z}/b\times O^\star_n\to(\mathbb{Z}/a)^\star$.
By Proposition \ref{AQ}, $\mbox{Aut}_{\tau_2}(O^\star_n)=\mbox{Aut}\,(O^\star_n)$, so
we achieve the factorization
$$\xymatrix{\mbox{Aut}\,({\mathbb Z}/a\rtimes_\tau(\mathbb{Z}/b\times O^\star_n))
\ar[dd]_{\psi}\ar[rr]^\eta&&\mbox{Aut}\,(H^{2k[\ell(\tau),2]}({\mathbb Z}/a\rtimes_\tau(\mathbb{Z}/b\times O^\star_n)) \\
\\
({\mathbb Z}/a)^\star\times\mbox{Aut}_{\tau_1}(\mathbb{Z}/b)\times\mbox{Aut}\,(O^\star_n)\ar[uurr]_{\eta'}}
$$ for all $k > 0$.
Because $$(\mathbb{Z}/(16\times 3^n))^\star/\pm\mbox{Aut}\,(O^\star_n)=
((\mathbb{Z}/16)^\star/\{\pm 1,\,\pm 9\})\times((\mathbb{Z}/3^n)^\star/\{\pm l^{k[\ell(\tau),2]};\;l\in(\mathbb{Z}/3^n)^\star\}),$$
we derive that
$O_{\mbox{\tiny Aut}\,(O^\star_n)}(16\times 3^n,2k[\ell(\tau),2])=2\times 3^{n-1}$.
 Then, the discussion above yields the main result.
\begin{thmX}\k Let $\gamma=(\gamma_1,\gamma_2) : {\mathbb Z}/b\times T^\star_n\to({\mathbb Z}/a)^\star$ and
$\tau=(\tau_1,\tau_2) : {\mathbb Z}/b\times O^\star_n\to({\mathbb Z}/a)^\star$ be actions
with $(a,b)=(ab,6)=1$ and $n\ge 3$, where $\gamma_1 : {\mathbb Z}/b\to({\mathbb Z}/a)^\star$,
$\gamma_2 : T^\star_n\to({\mathbb Z}/a)^\star$ and $\tau_1: \mathbb{Z}/b\to(\mathbb{Z}/a)^\star$,
$\tau_2 : O^\star_n\to(\mathbb{Z}/a)^\star$ are appropriate restrictions
of $\gamma$ and $\tau$, respectively. Then$:$

\vspace{3mm}

{\em (1)} $\mbox{\em card}\,{\mathcal K}^{2k[\ell(\gamma),2]-1}_{\mathbb{Z}/a\rtimes_\gamma
(\mathbb{Z}/b\times T^\star_n)}/_\simeq=2^{t+t'+1}3^{n_0}O(a,k[\ell(\gamma),2])
O_{\mbox{\tiny\em Aut}_{\gamma_1}(\mathbb{Z}/b)}(b,k[\ell(\gamma),2])O(3^{n-n_0},$

$k[\ell(\gamma),2])^{-1}$ for some $0\le t\le 2$ and $0\le t'\le 1$$;$

\vspace{3mm}

{\em (2)} $\mbox{\em card}\,{\mathcal K}^{2k[\ell(\tau),2]-1}_{\mathbb{Z}/a\rtimes_\tau(\mathbb{Z}/b\times O^\star_n)}/_\simeq=2^{t+1}\times 3^{n-1}
O(a,k[\ell(\tau),2])O_{\mbox{\tiny\em Aut}_{\tau_1}(\mathbb{Z}/b)}(b,k[\ell(\tau),2])$ for

some $0\le t\le 1$.
\end{thmX}
We point out that the numbers $t,t'$ above are given by \cite[Proposition 2.2]{GG2} and the orders
$O_{\mbox{\tiny Aut}_{\gamma_1}\,({\mathbb Z}/b)}(b,k[\ell(\gamma),2])$,
$O_{\mbox{\tiny Aut}_{\tau_1}\,
({\mathbb Z}/b)}(b,k[\ell(\tau),2])$ are determined by \cite[Corollary 2.3]{GG2}.

\bigskip

\par Now, let $\gamma=(\gamma_1,\gamma_2) : {\mathbb Z}/b\times T^\star_n\to
({\mathbb Z}/a)^\star$ and $\tau=(\tau_1,\tau_2) :
{\mathbb Z}/b\times O^\star_n\to ({\mathbb Z}/a)^\star$ be actions with
$\ell(\gamma),\ell(\tau)\le 2$. Then, of course
$2\ell(\gamma),2[\ell(\tau),2]\le 4$, a {\em fortiori} $\gamma_2$ is trivial
and the groups ${\mathbb Z}/a\rtimes_\gamma({\mathbb Z}/b\times T^\star_n)$, and ${\mathbb Z}/a\rtimes_\tau({\mathbb Z}/b\times O^\star_n)$
act on a $CW$-complex $X(4k-1)$ for any $k\ge 1$. Furthermore, as it was observed in \cite{GG2},
$\mbox{Aut}_{\gamma_1}\,({\mathbb Z}/b)=\mbox{Aut}_{\tau_1}\,({\mathbb Z}/b)=({\mathbb Z}/b)^\star$ .
Hence, ${\mathbb Z}/a\rtimes_\gamma({\mathbb Z}/b\times T^\star_n)\cong{\mathbb Z}/ab\times T^\star_n$
and ${\mathbb Z}/a\rtimes_\tau({\mathbb Z}/b\times O^\star_n)\cong{\mathbb Z}/ab\rtimes_{\tau'}O^\star_n$
with the action $\tau' : O^\star_n\stackrel{\tau_2}{\to}({\mathbb Z}/a)^\star\hookrightarrow({\mathbb Z}/ab)^\star$,
respectively. Then, in the light of \cite[Proposition 2.2]{GG2}, we are in a position to deduce
\begin{CorX}\k Let $\gamma : {\mathbb Z}/b\times T^\star_n\to({\mathbb Z}/a)^\star$ and $\tau :
{\mathbb Z}/b\times O^\star_n\to({\mathbb Z}/a)^\star$ be actions with $(a,b)=(ab,6)=1$
and $\ell(\gamma),\ell(\tau)\le 2$. Then$:$

\vspace{2mm}

{\em (1)}  $\mbox{\em card}\,{\mathcal K}^{4k-1}_{\mathbb{Z}/a\rtimes_\gamma(\mathbb{Z}/b\times T^\star_n)}=2\times 3^n\mbox{\em card}\,{\mathcal K}_{{\mathbb Z}/ab}^{4n-1}/_\simeq$

\vspace{1mm}

and

\vspace{1mm}

{\em (2)} $\mbox{\em card}\,{\mathcal K}^{4k-1}_{\mathbb{Z}/a\rtimes_\tau(\mathbb{Z}/b\times O^\star_n)}=2\times 3^{n-1}\mbox{\em card}\,{\mathcal K}_{{\mathbb Z}/ab}^{4n-1}/_\simeq.$
\end{CorX}
We point out that $\mbox{card}\,{\mathcal K}_{{\mathbb Z}/ab}^{4n-1}/_\simeq$ as the number of
homotopy types of $(4n-1)$-lens spaces has been fully described in \cite{GG}.

\renewcommand{\thesection}{{}}
\section{}
\renewcommand{\thesection}{\arabic{section}}

{\large\bf 3.\ Groups of self homotopy equivalences.} Let $\mu$ be a free and
cellular action of a finite group $G$ on a $CW$-complex $X(2k-1)$.
Write $\tilde{\eta} : \mbox{Aut}\,(G)\to({\mathbb Z}/|G|)^\star/\{\pm1\}$
for the composition of the anti-homomorphism $\eta : \mbox{Aut}\,(G)\to H^{2k}(G)=({\mathbb Z}/|G|)^\star$
 considered in the previous section with the quotient map $({\mathbb Z}/|G|)^\star\to({\mathbb Z}/|G|)^\star/\{\pm1\}$.
Then, by means of \cite[Proposition 3.1]{GG} (see also \cite[Theorem 1.4]{Sm}), the group
${\mathcal E}(X(2k-1)/\mu)$ of homotopy classes of self homotopy equivalences
for the space form $X(2k-1)/\mu$ is independent of the action $\mu$ of the
group $G$ as isomorphic to the kernel of the map
$\tilde{\eta} : \mbox{Aut}\,(G)\to({\mathbb Z}/|G|)^\star/\{\pm 1\}$ for all $n\ge 1$ provided $|G|>2$.
Whence, we simply write ${\mathcal E}(X(2k-1)/G)$ for this group.
\par Let $\alpha : G\to({\mathbb Z}/a)^\star$ be an action, $(a,|G|)=1$ and $2d$ a period of $G$.
Then, in virtue of \cite[Theorem 1.4]{Sm}, one gets $${\mathcal E}\,(X(2kd-1)/({\mathbb Z}/a\rtimes_\alpha G))\cong\left\{\begin{array}{ll}{\mathbb Z}/2,&\mbox{if}\;a|G|\le 2;\\
{\mathcal E}\,(X(2kd-1)/G)),&\mbox{if}\;a|G|> 2\;\mbox{and}\,a\le 2.\\
\end{array}\right.$$
By \cite{GG2}, the following generalization of \cite[Theorem 1.8]{Sm} holds.
\begin{proX}\k Let the group ${\mathbb Z}/a\rtimes_\alpha G$ with $(a,|G|)=1$ acts freely and
celullary on a $CW$-complex $X(2k[\ell(\alpha),d]-1)$ for $n\ge 1$, where $2d$ is a period of $G$.
Then, there are isomorphisms$:$
{\small
$${\mathcal E}\,(X(2k[\ell(\alpha),d]-1)/({\mathbb Z}/a\rtimes_\alpha G))\cong\left\{\begin{array}{ll}{\mathbb Z}/2,&\mbox{if}\;a|G|\le 2;\\
{\mathcal E}\,(X(2k[\ell(\alpha),d]-1)/G)),&\mbox{if}\,a|G|> 2\;\mbox{and}\;a\le 2;\\
\end{array}\right.$$ \label{PI}}
however for $a|G|>2$ with $a > 2$ it holds
$${\mathcal E}(X(2k[\ell(\alpha),d]-1)/({\mathbb Z}/a\rtimes_\alpha G))\cong$$

\noindent
{\small $\left\{\begin{array}{ll}\mbox{\em Der}_\alpha\,(G,\mathbb{Z}/a)
\rtimes({\mathcal E}\,(X(2k[\ell(\alpha),d]-1)/({\mathbb Z}/a))
\times{\mathcal E}_\alpha\,(X(2n[\ell(\alpha),d]-1)/G)),&\mbox{if}\;|G|>2;\\
{\mathbb Z}/a\rtimes{\mathcal E}\,(X(2k[\ell(\alpha),d]-1)/({\mathbb Z}/a)),&\mbox{if}\;|G|\le 2,\\
\end{array}\right.$}

\vspace{1mm}

\noindent
where ${\mathcal E}_\alpha\,(X(2n[\ell(\alpha),d]-1)/G))$ is the subgroup of ${\mathcal E}\,(X(2n[\ell(\alpha),d]-1)/G))$
determined by the subgroup $\mbox{\em Aut}_\alpha\,(G)\subseteq\mbox{\em Aut}\,(G)$.
\end{proX}

\bigskip

\par The paper \cite[Section 3]{GG1} deals with the group ${\mathcal E}\,(X(2k[\ell(\alpha),d]-1)/({\mathbb Z}/a))$,
however the group ${\mathcal E}_\alpha\,(X(2k[\ell(\alpha),d]-1)/G)$ consists of all automorphisms
$\varphi\in\mbox{Aut}_\alpha\,(G)$ with $\varphi^\ast=\pm\mbox{id}_{({\mathbb Z}/|G|)^\star}$ provided $|G|>2$.

\bigskip

\par If now $\gamma=(\gamma_1,\gamma_2) : \mathbb{Z}/b\times T^\star_n\to(\mathbb{Z}/a)^\star$ and
$\tau=(\tau_1,\tau_2) : \mathbb{Z}/b\times O^\star_n\to(\mathbb{Z}/a)^\star$ are actions considered in the previous section
then it holds $|\mathbb{Z}/b\times T^\star_n|=8\times 3^nb>2$ and $|\mathbb{Z}/b\times O^\star_n|=16\times 3^nb>2$
for the order of those groups. Furthermore,
$${\mathcal E}_\gamma\,(X(2k[\ell(\gamma),2]-1)/({\mathbb Z}/b)\times T^\star_n)\cong$$
$${\mathcal E}_{\gamma_1}\,(X(2k[\ell(\tau),2]-1)/({\mathbb Z}/b))\times{\mathcal E}_{\gamma_2}\,(X(2k[\ell(\gamma),2]-1)/T^\star_n)$$
and
$${\mathcal E}_\tau\,(X(2k[\ell(\tau),2]-1)/({\mathbb Z}/b)\times O^\star_n)\cong$$
$${\mathcal E}_{\tau_1}\,(X(2k[\ell(\tau),2]-1)/({\mathbb Z}/b))\times{\mathcal E}_{\tau_2}\,(X(2k[\ell(\tau),2]-1)/O^\star_n).$$
But, the group ${\mathcal E}_{\gamma_1}\,(X(2k[\ell(\gamma),2]-1)/({\mathbb Z}/b))$ and
${\mathcal E}_{\tau_1}\,(X(2k[\ell(\tau),2]-1)/({\mathbb Z}/b))$
has been fully described in \cite[Theorem 3.2]{GG2}.
\par To study the group ${\mathcal E}_{\gamma_2}\,(X(2k[\ell(\gamma),2]-1)/T^\star_n)$, we recall that by \cite{Sw} any
automorphism $\varphi\in\mbox{Aut}\,(Q_8)\cong S_4$ induces the identity
map on the cohomology group $H^{k[\ell(\gamma),2]}(Q_8)$ and by Proposition \ref{AQ},
$\mbox{Aut}_{\gamma_2}\,(T^\star_n)\cong S_4\times\mathbb{Z}/3^{n-n_0}$
provided $\mbox{Ker}\,\gamma_2=Q_8\rtimes_{\alpha_n}\mathbb{Z}/3^{n_0}$.
Then, we easily derive an isomorphism
$${\mathcal E}_{\gamma_2}\,(X(2k[\ell(\gamma),2]-1)/T^\star_n)\cong S_4\times\mathbb{Z}/\bigg(\frac{3^{n-n_0}}{(3^{n-n_0},k[\ell(\gamma),2])}\bigg).$$
\par Now, to move to the group ${\mathcal E}_{\tau_2}\,(X(2k[\ell(\tau),2]-1)/O^\star_n)$, we first recall that by Proposition \ref{AQ},
$\mbox{Aut}_{\tau_2}(O^\star_n)=\mbox{Aut}\,(O^\star_n)$. Because of the isomorphism
$H^{2k[\ell(\gamma),2]}(O_n^\star)\cong\mathbb{Z}/(16\times 3^nab)$ from Section 1, we must study all automorphisms
$\varphi\in\mbox{Aut}\,(O^\star_n)$ with $\varphi^\ast=\mp\mbox{id}_{\mathbb{Z}/(16\times 3^nab)}$.


\par Let $\mbox{Aut}^0\,(O^\star_n)=\{\varphi\in\mbox{Aut}\,(O^\star_n);\,\varphi^\ast=\mbox{id}_{\mathbb{Z}/(16\times 3^nab)}\}$.
Consider the automorphisms $\varphi,\psi\in\mbox{Aut}\,(O^\star_n)$ defined on
generators (according to the presentation of $O^\star_n$ given in Section 1) by:
$\varphi(P)=P$, $\varphi(Q)=Q$, $\varphi(X)=X$, $\varphi(R)=R^{-1}$ and $\psi(P)=P$,
$\psi(Q)=Q$, $\psi(X)=X^4$, $\varphi(R)=R$ and write $\big<\varphi,\psi\big>$ for the subgroup
of $\mbox{Aut}\,(O^\star_n)$ generated by $\varphi$ and $\psi$. It is easy to
check that there is an isomorphism $\big<\varphi,\psi\big>\cong\mathbb{Z}/2\oplus\mathbb{Z}/3^{n-1}$
and $\big<\varphi,\psi\big>\cap\mbox{Inn}\,(O^\star_n)=E$, the trivial subgroup of $O^\star_n$.
Then, by the results of \cite{GG3} and order arguments, there is the splitting short exact sequence
$$1\to\mbox{Inn}\,(O^\star_n)\longrightarrow\mbox{Aut}(O_n^\star)
\longrightarrow\mathbb{Z}/2\oplus\mathbb{Z}/3^{n-1}\to 1.$$
Consequently, a simple calculation, by means of the list of elements in $\mbox{Aut}\,(O^\star_n)$
presented in \cite{GG3} and considerations in the first paragraph on page 14 provides an isomorphism
$$\mbox{Aut}^0(O^\star_n)\cong\mbox{Inn}\,(O^\star_n)\rtimes\mathbb{Z}/\bigg(\frac{3^{n-1}}{(3^{n-1},k[\ell(\tau),2])}\bigg).$$
Because the restriction of any automorphism of $\mbox{Aut}\,(O^\star_n)$ to the subgroup $Q_8$
induces the identity in cohomology at dimension multiple of $4$, there is no element of $\mbox{Aut}\,(O^\star_n)$
which induces the minus identity in cohomology at dimension $2k[\ell(\tau),2]$.
Since $\mathcal{Z}(O^\star_n)=\mathbb{Z}/2$ and so $\mbox{Inn}\,(O^\star_n)\cong O^\star_n/(\mathbb{Z}/2)=O_n$
(the group considered in \cite{GG3}) and consequently, we derive an isomorphism
$${\mathcal E}_{\tau_2}\,(X(2k[\ell(\tau),2]-1)/O^\star_n)
\cong O_n\rtimes\mathbb{Z}/\bigg(\frac{3^{n-1}}{(3^{n-1},k[\ell(\tau),2])}\bigg).$$
\par Finally, by Proposition \ref{PI} and the consideration above, we can close the paper with
\begin{thmX}\k Let $\gamma=(\gamma_1,\gamma_2) : \mathbb{Z}/b\times T^\star_n\to(\mathbb{Z}/a)^\star$ $($resp.\
$\tau=(\tau_1,\tau_2) : \mathbb{Z}/b\times O^\star_n\to(\mathbb{Z}/a)^\star$$)$ be an action with
$(a,b)=(ab,6)=1$ for $n\ge 3$.
If the group ${\mathbb Z}/a\rtimes_\gamma({\mathbb Z}/b\times
T^\star_n)$ $($resp.\ ${\mathbb Z}/a\rtimes_\tau({\mathbb Z}/b\times O^\star_n)$$)$
acts freely and cellularly on a $CW$-complex $X(2k[\ell(\gamma),2]-1)$
$($resp.\ $X(2k[\ell(\tau),2]-1)$$)$ then
\begin{multline*}{\mathcal E}(X(2k[\ell(\gamma),2]-1)/({\mathbb Z}/a\rtimes_\gamma({\mathbb Z}/b\times
T^\star_n)))\cong\mbox{\em Der}_\gamma\,
({\mathbb Z}/b\times T^\star_n,{\mathbb Z}/a)\rtimes\\({\mathcal E}(X(2k[\ell(\gamma_1),2]-1)/({\mathbb Z}/a))\times
{\mathcal E}_{\gamma_1}\,(X(2k[\ell(\beta),2]-1)/({\mathbb Z}/b))\times
S_4\times\\\mathbb{Z}/\bigg(\frac{3^{n-n_0}}{(3^{n-n_0},k[\ell(\gamma),2])}\bigg)\\\end{multline*}

\vspace{-2cm}

\begin{multline*}(\mbox{resp.}\;\;{\mathcal E}(X(2k[\ell(\tau),2]-1)/({\mathbb Z}/a\rtimes_\tau({\mathbb Z}/b\times O^\star_n)))\cong\mbox{\em Der}_\tau\,
({\mathbb Z}/b\times O^\star_n,{\mathbb Z}/a)\rtimes\\({\mathcal E}(X(2k[\ell(\tau_1),2]-1)/({\mathbb Z}/a))\times
{\mathcal E}_{\tau_1}\,(X(2k[\ell(\beta),2]-1)/({\mathbb Z}/b))\times
O_n\rtimes\\\mathbb{Z}/\bigg(\frac{3^{n-1}}{(3^{n-1},k[\ell(\tau),2])}\bigg)).\\\end{multline*}
\end{thmX}

Thus, in the light of Proposition \ref{D}, the groups ${\mathcal E}(X(2k[\ell(\gamma),2]-1)/({\mathbb Z}/a\rtimes_\gamma({\mathbb Z}/b\times
T^\star_n)))$ and ${\mathcal E}(X(2k[\ell(\tau),2]-1)/({\mathbb Z}/a\rtimes_\tau({\mathbb Z}/b\times O^\star_n)))$ have been fully described.

\newpage

\vspace*{-1.5cm}

\renewcommand{\baselinestretch}{.07}

\large

\normalsize

{}

\vspace{7mm}

\noindent{\bf Faculty of Mathematics and Computer Science\\
Nicolaus Copernicus University\\
Chopina 12/18, 87-100 Toru\'n, Poland}\\
e-mail: marek@mat.uni.torun.pl

\vspace{5mm}

\noindent{\bf Department of Mathematics-IME\\
University of S\~ao Paulo\\
Caixa Postal 66.281-AG.\ Cidade de S\~ao Paulo\\
05311-970 S\~ao Paulo, Brasil}\\
e-mail: dlgoncal@ime.usp.br


\begin{thebibliography}{99}
\bibitem{AM} A.\ Adem, R.J.\ Milgam, {\em Cohomology of Finite Groups}, Springer-Verlag,
New York-Heidelberg-Berlin (1994).
\bibitem{AD} ----------------- , J.F.\ Davis, {\em Topics in Transformation Groups}, Handbook of geometric
topology, 1-54, North-Holland, Amsterdam (2002).
\bibitem{CE} H.\ Cartan, S.\ Eilenberg, {\em Homological Algebra}, Princeton, New Jersey (1956).
\bibitem{GG} M.\ Golasi\'nski, D.L.\ Gon\c{c}alves, {\em Homotopy spherical space forms -
a numerical bound for homotopy types}, Hiroshima Math.\ J.\ {\bf
31} (2001), 107-116.
\bibitem{GG1} ---------------- , {\em Spherical space forms -homotopy types and self-equivalences}, 
Progr.\ Math., {\bf 215} Birkh\"auser, Basel (2004), 153-165.
\bibitem{GG2} ---------------- , {\em Spherical space forms - homotopy types and self-equivalences for the groups $\mathbb{Z}/a\rtimes\mathbb{Z}/b$ and $\mathbb{Z}/a\rtimes(\mathbb{Z}/b\times\mathbb{Q}_{2^i})$},
Topology and its Appl.\ {\bf 146-147} (2005), 451-470.
\bibitem{GG3} ---------------- , {\em Automorphisms of generalized $($binary$)$ polyhedral groups}, (submitted).
\bibitem{GG4} ---------------- , {\em Spherical space forms - homotopy types and self-equivalences for
the groups $(\mathbb{Z}/a\rtimes\mathbb{Z}/b)\times
SL_2(\mathbb{F}_p)$}, (submitted).
\bibitem{R} J.W.\ Rutter, {\em Spaces of homotopy self-equivalences. A survey}, Lecture Notes in Math.\ {\bf 1662}, Springer-Verlag, Berlin (1997).
\bibitem{Sm} D.\ Smallen, {\em The group of self-equivalences of certain complexes}, Pacific J.\ Math.\ {\bf 54} (1974), 269-276.
\bibitem{Sw} R.G.\ Swan, {\em Periodic resolutions for finite groups}, Ann.\ of Math.\ {\bf 72}, no.\ 2 (1960), 267-291.
\bibitem{Th1} C.B.\ Thomas, {\em The oriented homotopy type of compact $3$-manifolds},
Proc.\ London Math.\ Soc.\ (3), {\bf 19} (1969), 31-44.
\bibitem{Wo} J.A.\ Wolf, {\em Spaces of constant curvature}, Publisher or Perish, Inc.\ Wilmington, Delaware, U.S.A.\ (1984).
\end{thebibliography}
\end{document}